\documentclass{amsart}
\usepackage[all]{xy}
\usepackage{lipsum}
\usepackage{amsfonts}
\usepackage{graphicx}
\usepackage{amsmath}
\usepackage{amssymb}
\usepackage[english]{babel}
\usepackage{tikz-cd}
\usepackage{hyperref}

\newtheorem{proposition}{Proposition}
\newtheorem{corollary}{Corollary}

\begin{document}

\title{Matricial Closure}

\author{Frank Murphy-Hernandez}
\address{Facultad de Ciencias, UNAM, Mexico City}
\email{murphy@ciencias.unam.mx}

\author{Francisco Raggi}
\address{Instituto de Matemáticas, UNAM, Mexico City}
\email{fraggi@matem.unam.mx}

\author{José Ríos}
\address{Instituto de Matemáticas, UNAM, Mexico City}
\email{ jrios@matem.unam.mx }


\date{\today}


\begin{abstract}
We study  a closure operator derived from the matrix endofunctor on the category of rings with unity. We investigate the invariance of various ring-theoretic properties under this operator. A key finding is the decisive nature of this operator: a property is either preserved for all rings or fails for all rings. This work provides a comprehensive analysis of the behavior of ring properties under matricial closure, including chain conditions, radical-theoretic properties, and other structural characteristics.
\end{abstract}

\maketitle

\section{Introduction}

The concept of a \emph{matricial closure} arises naturally from the theory of ultramatricial algebras. By considering a very specific type of direct system—one indexed by the powers of a fixed natural number $n$—we are able to define a categorical closure operator. This construction not only applies to algebras over a field but, more generally, to arbitrary rings, thus extending the concept beyond its traditional setting.

We prove that the matricial closure acts as a closure operator on the category of rings with unity, where the category is viewed as a large preorder induced by monomorphisms. Although ultramatricial algebras have been extensively studied in the context of the $K_0$ functor and ordered abelian groups with order units, our focus lies elsewhere.

First, we define the matricial closure not merely as an operator but as an endofunctor. We then show that it satisfies closure properties: an inflationary natural transformation that is monic, and an idempotent property realized through a natural isomorphism. Furthermore, we demonstrate that this closure commutes with certain categorical constructions, such as forming group rings and polynomial rings.

We also examine two naturally arising functors between the corresponding categories of left modules: one induced by the inflationary ring morphism, and another defined as a direct limit of Morita equivalences. Although neither functor is generally an equivalence, the latter preserves certain properties, such as simplicity and faithfulness.

Additionally, we analyze the lattice of ideals of the matricial closure, establishing a lattice isomorphism between the two-sided ideals of the original ring and those of its closure—analogous to the classical matrix ring case—while also showing that the lattice of left ideals can become significantly larger.

Finally, we investigate various chain and dimension conditions, showing that the matricial closure typically fails to satisfy them, which often makes it a pathological ring. We conclude by studying its behavior concerning the Jacobson radical, prime and semiprime rings, V-rings, von Neumann regular rings, and the invariant basis number property.

\section{Preliminaries}

Throughout this paper all rings will be assumed associative with 1. For a ring $R$ and positive integer $n$, $M_n(R)$ denotes the ring of matrices of $n \times n$ with coefficients in $R$, the 1 of this ring will be denoted by $I_n$ and satisfies $(I_n)_{ij} = \delta_{ij}$ for any $1 \leq i, j \leq n$ where $\delta$ is the Kronecker’s delta. The category of all rings with morphism that preserve the 1 is denoted by $\mathfrak{R}$. All $R$-modules will be left unitary modules over $R$.

Our reference for dimensions are [5] and [6]. For $R$ a ring and $M$ a left $R$-module a set $\{M_i\}_{i \in I}$ of submodules of $M$ are called independent, if $M_j \cap \sum_{i \in I, i \neq j} M_i = 0$ for all $j \in I$. If there is a maximal finite independent set of submodules of $M$, the cardinality of this set is the uniform dimension of $M$ and is well defined, in this case we will say that $M$ has finite uniform dimension, we also say that $R$ has left finite uniform dimension if it has finite uniform dimension as left module over itself.

For $R$ a ring and $M$ a left $R$-module, we define the Krull dimension of $M$, denoted by $d_K(M)$, recursively $d_K(0) = -\infty$, $d_K(M) = 0$ if $M$ is artinian, and $d_K(M) = \alpha$ if $d_K(M) \neq \beta$ for any ordinal $\beta < \alpha$ and in any descending chain of submodules of $M$ all but a finite number of factors have Krull dimension less than $\alpha$. A ring $R$ has left Krull dimension if it has Krull dimension as left module over itself. A classic result in the theory is that: Let $R$ be a ring and let $M$ be a left $R$-module. If $M$ has Krull dimension then $M$ has finite uniform dimension. As its corollary: Let $R$ be a ring. If $R$ has Krull dimension then $R$ has finite uniform dimension.

Let $R$ be a ring and let $M$ be a left $R$-module. Recall that $M$ has projective dimension $n$, denoted by $\text{pd}(M) = n$, if there is an exact sequence $0 \to P_n \to \dots \to P_0 \to M \to 0$ with $P_i$ projective for $i = 0, \dots, n$ and there is no such as that for any $0 \leq k < n$, immediately $M$ is projective if and only if $\text{pd}(M) = 0$.

\subsection{Dimensions in Ring Theory}

Our reference for dimensions are [5] and [6]. For $R$ a ring and $M$ a left $R$-module a set $\{M_i\}_{i \in I}$ of submodules of $M$ are called independent, if $M_j \cap \sum_{i \in I, i \neq j} M_i = 0$ for all $j \in I$. If there is a maximal finite independent set of submodules of $M$, the cardinality of this set is the uniform dimension of $M$ and is well defined, in this case we will say that $M$ has finite uniform dimension, we also say that $R$ has left finite uniform dimension if it has finite uniform dimension as left module over itself.

For $R$ a ring and $M$ a left $R$-module, we define the Krull dimension of $M$, denoted by $d_K(M)$, recursively $d_K(0) = -\infty$, $d_K(M) = 0$ if $M$ is artinian, and $d_K(M) = \alpha$ if $d_K(M) \neq \beta$ for any ordinal $\beta < \alpha$ and in any descending chain of submodules of $M$ all but a finite number of factors have Krull dimension less than $\alpha$. A ring $R$ has left Krull dimension if it has Krull dimension as left module over itself. A classic result in the theory is that: Let $R$ be a ring and let $M$ be a left $R$-module. If $M$ has Krull dimension then $M$ has finite uniform dimension. As its corollary: Let $R$ be a ring. If $R$ has Krull dimension then $R$ has finite uniform dimension.

Let $R$ be a ring and let $M$ be a left $R$-module. Recall that $M$ has projective dimension $n$, denoted by $\text{pd}(M) = n$, if there is an exact sequence $0 \to P_n \to \dots \to P_0 \to M \to 0$ with $P_i$ projective for $i = 0, \dots, n$ and there is no such as that for any $0 \leq k < n$, immediately $M$ is projective if and only if $\text{pd}(M) = 0$.

\subsection{K-Theory}

Our reference for K-theory are [1] and [7]. A ring $R$ is said to have Invariant Basis Number, if for any pair of naturals $n$ and $m$, $R^n \cong R^m$ as left $R$-modules implies $n = m$. A ring is called left projective free, if all finite generated projective modules are free.

We define $\mathcal{P}(R)$ as the class of all finitely generated projective left $R$-modules for any ring $R$, and so we may consider the monoid $\mathcal{M}(R)$ of the isomorphism classes of $\mathcal{P}(R)$ and the coproduct as operation, then an element $[P]$ in $\mathcal{M}(R)$ is the equivalence class of finitely generated projective left $R$-modules and therefore $[P] = [P']$ if and only if $P \cong P'$ as left $R$-modules. The operation is given by $[P] + [Q] = [P \oplus Q]$. It is clear that in this way $\mathcal{M}(R)$ is a commutative monoid. Finally we define by $K_0(R)$ the Grothendieck group of $\mathcal{M}(R)$. The functor $K_0$ relates with the invariant basis number in the following way: Let $R$ be a ring and let $\lambda : \mathbb{Z} \to K_0(R)$ be the group morphism induced by $\lambda(1) = [R]$. Then $R$ has the invariant basis number if and only if $\lambda$ is a monomorphism. Moreover $\lambda$ restricted to the submonoid generated by $[R]$ to $\mathbb{N}$ is a monoid isomorphism if and only if $R$ is left projective free. With respect to direct limits, the functor $K_0$ commutes with direct limits, that is, let $(R_i, \alpha_i^j : R_i \to R_j)_{i \leq j \in I}$ be a direct system of rings. Then $\varinjlim K_0(R_i) \cong K_0(\varinjlim R_i)$.

\subsection{Classes of rings}

For general aspects of ring theory our references are [4] and [8]. Now we give certain properties that the rings may have and we will test if are preserved by our endofunctor, also we will give equivalences that will make the proof in the easiest way. Let $R$ be a ring. Recall that the Jacobson radical of $R$, denoted by $J(R)$, is the intersection of all maximal left ideals of $R$. A $x \in R$ is left quasiregular in $R$ if $1 - x$ is left invertible in $R$, between the many characterizations of $J(R)$ is that $J(R)$ is the unique left ideal which all its elements are left quasiregular and is maximal with respect to this property. A well known result is the for any $n$ positive integer, $J(M_n(R)) = M_n(J(R))$. A ring $R$ such that $J(R) = 0$ is called semisimple, and a semisimple artinian ring is such that is left artinian and semisimple.

A ring $R$ is called left perfect if any left $R$-module has a projective cover. The equivalence that we will use is that a ring $R$ is left perfect if and only if $R$ is a left max ring (all non zero left $R$-modules have a maximal submodule) and $R/J(R)$ is semisimple artinian. A ring $R$ is called semiperfect, if any left simple $R$-module has a projective cover, and this is equivalent to $R/J(R)$ to be artinian semisimple and idempotents lift modulo $J(R)$.

A ring $R$ is called von Neumann regular ( see [2]) if for any $a \in R$ exists $x \in R$ such that $a = axa$. A classic result is that if $R$ is a von Neumann regular ring, then $M_n(R)$ is von Neumann regular ring for any natural $n$. A ring $R$ is called left V-ring if all simple left $R$-modules are injective. A ring $R$ is called quasi Frobenius if it is left noetherian and left self injective.

A ring $R$ is called left primitive ring if it has a simple faithful left $R$-module. An equivalence is that a ring $R$ is left primitive if and only if $R$ is isomorphic to a dense subring of the ring of endomorphisms of a left vector space over a division ring $D$. A left full linear ring $R$ is the ring of all linear transformations of an infinite dimensional left vector space over a division ring $D$, this is equivalent to ask to $R$ to be von Neumann regular, left self-injective with non zero left socle.

Let $R$ be a ring and let $I$ be an ideal of $R$, $I$ is semiprime if for any $a, x \in R$, $axa \in I$ implies $a \in I$. A ring is called semiprime if the zero ideal is semiprime. Let $R$ be a ring and let $I$ be an ideal of $R$, $I$ is prime if for any $a, b, x \in R$, $axb \in I$ implies $a \in I$ or $b \in I$. A ring is called prime if the zero ideal is prime. Trivially all prime rings are semiprime rings.

\section{Matricial Closure}

\subsection{Definition}

Let $n$ be a positive integer, we define the assignation $M_n : \mathfrak{R} \to \mathfrak{R}$ given by: for any ring $R$ we assign it $M_n(R)$ and for any ring morphism $f : R \to S$, $M_n(f) : M_n(R) \to M_n(S)$ is given by: for any $A \in M_n(R)$, $M_n(f)(A)_{ij} = f(A_{ij})$ with $1 \leq i, j \leq n$.

\begin{proposition}
Let $n$ be a positive integer. Then $M_n : \mathfrak{R} \to \mathfrak{R}$ is a functor.
\end{proposition}

\textbf{Proof.} Obviously $M_n$ sends rings into rings, what is left is to verify that sends morphisms into morphisms. Let $f : R \to S$ be a ring morphism, $A, B \in M_n(R)$ and $1 \leq i, j \leq n$, first $M_n(f)(A+B)_{ij} = f((A+B)_{ij}) = f(A_{ij} + B_{ij}) = f(A_{ij}) + f(B_{ij}) = M_n(f)(A)_{ij} + M_n(f)(B)_{ij}$ which means that $M_n(f)(A+B) = M_n(f)(A) + M_n(f)(B)$. Second $M_n(f)(AB)_{ij} = f((AB)_{ij}) = f(\sum_{k=1}^n A_{ik}B_{kj}) = \sum_{k=1}^n f(A_{ik})f(B_{kj}) = \sum_{k=1}^n M_n(f)(A)_{ik}M_n(f)(B)_{kj} = (M_n(f)(A)M_n(f)(B))_{ij}$ which means that $M_n(f)(AB) = M_n(f)(A)M_n(f)(B)$. And finally $M_n(f)(I_n)_{ij} = f((I_n)_{ij}) = f(\delta_{ij}) = \delta_{ij} = (I_n)_{ij}$.

And for any positive integer $n$ we have also a natural transformation $\eta^n : 1_{\mathfrak{R}} \to M_n$ given: for any $R$ ring and $a \in R$, $\eta^n_R : R \to M_n(R)$ is given by $\eta^n_R(a)_{ij} = a \delta_{ij}$ with $1 \leq i, j \leq n$.

\begin{proposition}
Let $n$ be a positive integer. Then $\eta^n : 1_{\mathfrak{R}} \to M_n$ is a natural transformation.
\end{proposition}

\textbf{Proof.} Let $f : R \to S$ be a ring morphism and $a \in R$, first $((\eta^n_S \circ f)(a))_{ij} = (\eta^n_S(f(a)))_{ij} = f(a)\delta_{ij}$ with $1 \leq i, j \leq n$, for the other side $((M_n(f) \circ \eta^n_R)(a))_{ij} = (M_n(f)(\eta^n_R(a)))_{ij} = f(\eta^n_R(a)_{ij}) = f(a \delta_{ij}) = f(a)\delta_{ij}$, so $M_n(f)\eta^n_R = \eta^n_S f$ and the diagram commutes
\[
\xymatrix{
R \ar[r]^f \ar[d]_{\eta^n_R} & S \ar[d]^{\eta^n_S} \\
M_n(R) \ar[r]_{M_n(f)} & M_n(S)
}
\]

It is well known that for any pair of positive integers $n$ and $m$, there is a natural isomorphism from $M_n \circ M_m$ to $M_{nm}$. So we fix $n$ and for any naturals $m, k$ consider $\alpha[R,n]^k_m = \bigcirc_{i=m}^{k-1} \eta^n_{M_{n^i}(R)}$ if $m < k$ and $\alpha[R,n]^k_m = 1_{M_{n^k}(R)}$ if $m = k$. The last composition was made identifying $M_k(M_m(R))$ with $M_{kn}(R)$ by the natural isomorphism mentioned before. For any ring $R$ and any positive integer $n$ we get a directed system over the naturals, we denote its direct limit by $MC_n(R)$. We will denote by $i[R]^n_k : M_{n^k}(R) \to MC_n(R)$ the canonical morphism in the direct limit, also note that this morphism is injective, then a ring monomorphism. Now consider a ring morphism $f : R \to S$ then we take the family of ring morphisms $\{M_{n^k}(f) : M_{n^k}(R) \to M_{n^k}(S)\}_{k \in \mathbb{N}}$ induced for the family of functors $\{M_{n^k}\}_{k \in \mathbb{N}}$ and composed with the family of morphism $\{i[S]^n_k\}_{k \in \mathbb{N}}$ so we get the following family of ring morphisms $\{i[S]^n_k \circ M_{n^k}(f) : M_{n^k}(R) \to MC_n(S)\}_{k \in \mathbb{N}}$, finally this family is compatible so by the universal property of the direct limit there exist a unique ring morphism $MC_n(f) : MC_n(R) \to MC_n(S)$.

\begin{proposition}
Let $n$ be a positive integer. Then $MC_n$ is a functor.
\end{proposition}

\textbf{Proof.} Again it is obvious that $MC_n : \mathfrak{R} \to \mathfrak{R}$ sends rings into rings. So we will focus in the morphisms, first notice that $\alpha[\_, n]^k_m : M_{n^m} \to M_{n^k}$ is natural transformation since it is a composition of many natural transformations. So for any ring morphism $f : R \to S$ we get the next commutative diagram
\[
\xymatrix{
M_{n^m}(R) \ar[r]^{M_{n^m}(f)} \ar[d]_{\alpha[R,n]^k_m} & M_{n^m}(S) \ar[d]^{\alpha[S,n]^k_m} \\
M_{n^k}(R) \ar[r]_{M_{n^k}(f)} & M_{n^k}(S)
}
\]
That is $M_{n^k}(f) \circ \alpha[S,n]^k_m = \alpha[R,n]^k_m \circ M_{n^m}(f)$. Now consider that the canonical injections satisfy $i[S]^n_m = i[S]^n_k \circ \alpha[S,n]^k_m$, so composing the last equality with $i[S]^n_k$ by the left we get $i[S]^n_k \circ M_{n^k}(f) \circ \alpha[S,n]^k_m = i[S]^n_k \circ \alpha[R,n]^k_m \circ M_{n^m}(f) = i[S]^n_m \circ M_{n^m}(f)$ then compatibility condition is satisfied so we get $MC_n(f)$. $\blacksquare$

We call $MC_n(R)$ the $n$-matricial closure of the ring $R$.

\subsection{Closure Properties}

\begin{proposition}[Inflatory]
Let $n, k$ be positive integers. Then $i^n_k : M_{n^k} \to MC_n$ is an injective natural transformation.
\end{proposition}

\textbf{Proof.} We wish that the following diagram to commute
\[
\xymatrix{
M_{n^k}(R) \ar[r]^{M_{n^k}(f)} & M_{n^k}(S) \\
MC_n(R) \ar[r]_{MC_n(f)} & MC_n(S)
}
\]
for any ring morphism $f : R \to S$, so we recall that $MC_n(f) \circ i[R]^n_k = i[S]^n_k \circ M_{n^k}(f)$ this by the universal property of the direct limit.

In particular we will denote by $i_n : 1_{\mathfrak{R}} \to MC_n$ the natural transformation above in the case that $k = 0$.

\begin{proposition}[Monotone]
Let $n$ be a positive integer, let $R$ and $S$ be rings and let $f : R \to S$ a ring monomorphism. Then $MC_n(f)$ is a ring monomorphism.
\end{proposition}

\textbf{Proof.} As the functor $M_{n^k}$ sends monomorphism into monomorphisms, we get proposition.

\begin{proposition}[Idempotent]
Let $n$ be a positive integer. Then there is a natural isomorphism between $MC_n$ and $MC_n^2$.
\end{proposition}

\textbf{Proof.} First we will prove that there is a natural isomorphism from $M_n \circ MC_n$ to $MC_n$. To do this, we build $\epsilon[n]_R : M_n(MC_n(R)) \to MC_n(R)$, as follows, first take an element $A \in M_n(MC_n(R))$, so for all $1 \leq i, j \leq n$, $A_{ij} = i[R]^n_{k_{ij}}(B_{ij})$ for some $k_{ij} \in \mathbb{N}$ and $B_{ij} \in M_{n^{k_{ij}}}(R)$ and we put $m = \max\{k_{ij} \mid 1 \leq i, j \leq n\}$. Without loss of generality we may assume that $k_{ij} = m$ for all $1 \leq i, j \leq n$, that is because we could consider that $\alpha[R,n]^{m}_{k_{ij}}(B_{ij}) \in M_{n^m}(R)$. Then we construct a new matrix $\widehat{A} \in M_n(M_{n^m}(R)) = M_{n^{m+1}}(R)$, given by $\widehat{A}_{ij} = B_{ij}$ for any $1 \leq i, j \leq n$, and finally we could define $\epsilon[n]_R(A) = i[R]^{m+1}_n(\widehat{A})$. Consider that $\epsilon_R(A) = 0$, in such case with the above notation we get that $i[R]^{m+1}_n(\widehat{A}) = 0$, but $i[R]^{m+1}_n$ is injective so $\widehat{A} = 0$ and this implies that $A = 0$, so we get the injectivity of $\epsilon[n]_R$. Now consider $B \in MC_n(R)$ with $B = i[R]^k_n(C)$ and $A \in M_n(MC_n(R))$ given by $A_{ij} = B \delta_{ij}$, thus $\widehat{A}_{ij} = C \delta_{ij}$ and so $\epsilon_R(A) = B$.

For the naturality, we take a ring morphism $f : R \to S$, so we would like this diagram to commute
\[
\xymatrix{
M_nMC_n(R) \ar[r]^{M_nMC_n(f)} \ar[d]_{\epsilon[n]_R} & M_nMC_n(S) \ar[d]^{\epsilon[n]_S} \\
MC_n(R) \ar[r]_{MC_n(f)} & MC_n(S)
}
\]
For this, we note that for an element $A \in M_nMC_n(R)$ and $1 \leq i, j \leq n$, $\widehat{M_nMC_n(f)(A)}_{ij} = M_{n^m}(f)(\widehat{A}_{ij}) = M_{n^{m+1}}(f)(\widehat{A})_{ij}$ which means that $M_n\widehat{MC_n(f)(A)} = M_{n^{m+1}}(f)(\widehat{A})$, and so
\begin{align*}
\epsilon[n]_S(M_nMC_n(f)(A)) &= i[S]^{m+1}_n(M_n\widehat{MC_n(f)(A)}) \\
&= i[S]^{m+1}_n(M_{n^{m+1}}(f)(\widehat{A})) \\
&= MC_n(f) i[R]^{m+1}_n(\widehat{A}) \\
&= MC_n(f) \epsilon[n]_R(A)
\end{align*}
that is the diagram commutativity.

Now we get $M_n(\epsilon[n]) : M_{n^2}MC_n \to M_nMC_n$ which is a natural isomorphism, so if we consider $\epsilon[n] \circ M_n(\epsilon[n]) : M_{n^2}MC_n \to MC_n$ we get another natural isomorphism, in fact, we construct a family of natural isomorphisms in the next fashion, recursively $\varepsilon[n,1] = \epsilon[n] : M_nMC_n \to MC_n$ and $\varepsilon[n,k+1] = \epsilon[n] \circ M_n(\varepsilon[n,k]) : M_{n^{k+1}}MC_n \to MC_n$. So we get a new natural isomorphism $\epsilon[n] : MC_nMC_n \to MC_n$ where $\varepsilon$ is the direct limit of the family of natural isomorphisms described.

As it proved above the matricial closure operator behaves as an operator closure in a partial ordered, but in this case we have consider the category $\mathfrak{R}$ with the preorder induced by the monomorphisms.
\section{Functorial Properties with respect to certain Constructions}

\begin{proposition}
Let $n$ be a positive integer and let $I$ be a set. Then there is a natural isomorphism $\upsilon[n,I] = \{\upsilon[n,I]_{(R_i), i \in I} : MC_n(\prod_{i \in I} R_i) \to \prod_{i \in I} MC_n(R_i)\}_{(R_i) \in \mathfrak{R}^I}$
\end{proposition}

\textbf{Proof.} Let $(R_i)_{i \in I}$ be a family of rings indexed over $I$, we construct $v[n,I]^1_{(R_i)_{i \in I}}$ from $M_n(\prod_{i \in I} R_i)$ to $\prod_{i \in I} M_n(R_i)$ given by $v[n]^1(A)(i)_{jk} = A(i)_{jk}$ for any $A \in M_n(\prod_{i \in I} R_i)$, $i \in I$ and $1 \leq j, k \leq n$, under the assumption that $v[n,I]^1$ is a natural isomorphism we may construct the desired natural isomorphism, so we take $v[n,I]^{m+1}_{(R_i)_{i \in I}} = v[n,I]^1_{(M_{n^m}(R_i))_{i \in I}} \circ M_n(v[n,I]^m_{(R_i)_{i \in I}})$ and as before we claim that the natural isomorphism $v[n,I]$ we looked for is the direct limit of the family of natural isomorphisms $\{v[n,I]^k\}_{k \in \mathbb{N}}$.

\begin{proposition}
Let $n$ be a positive integer. Then there is a natural isomorphism $\iota[n] = \{\iota[n]_R : MC_n(R)[x] \to MC_n(R[x])\}_{R \in \mathfrak{R}}$.
\end{proposition}

\textbf{Proof.} Let $R$ be a ring and we contract $\iota[n]^1_R : M_n(R)[x] \to M_n(R[x])$ given by $\theta[n]^1_R(A x)_{ij} = A_{ij} x$ for any $A \in M_n(R)$, $1 \leq i, j \leq n$, again under the assumption that $\iota[n]^1$ is a natural isomorphism we build the desired natural isomorphism, then $\iota[n]^{m+1}_R = \iota[n]^1_{M_n(R)} \circ M_n(\iota[n]^m_R)$. And again we get the natural isomorphism $\iota[n]$ as the direct limit of the family of natural isomorphism $\{\iota[n]^k\}_{k \in \mathbb{N}}$.

\begin{proposition}
Let $n$ be a positive integer and let $G$ be a group. Then there is a natural isomorphism $\theta[n] = \{\theta[n]_R : MC_n(R)[G] \to MC_n(R[G])\}_{R \in \mathfrak{R}}$.
\end{proposition}

\textbf{Proof.} Let $R$ be a ring and we contract $\theta[n]^1_R : M_n(R)[G] \to M_n(R[G])$ given by $\theta[n]^1_R(A g)_{ij} = A_{ij} g$ for any $A \in M_n(R)$, $1 \leq i, j \leq n$ and $g \in G$, again under the assumption that $\theta[n]^1$ is a natural isomorphism we build the desired natural isomorphism, then $\theta[n]^{m+1}_R = \theta[n]^1_{M_n(R)} \circ M_n(\theta[n]^m_R)$ and again we get the natural isomorphism $\theta[n]$ as the direct limit of the family of natural isomorphism $\{\theta[n]^k\}_{k \in \mathbb{N}}$.

\begin{corollary}
Let $n$ be a positive integer. Then there is a natural isomorphism $\Theta[n] = \{\Theta[n]_R : MC_n(R)[x, x^{-1}] \to MC_n(R[x, x^{-1}])\}_{R \in \mathfrak{R}}$.
\end{corollary}

\begin{proposition}
Let $n$ be a positive integer and let $R$ be a ring. Then $Cen(MC_n(R)) \cong Cen(R)$.
\end{proposition}

\textbf{Proof.} Consider that the next ring morphism $f : Cen(R) \to Cen(MC_n(R))$ where $f$ is $i_n$ restricted to $Cen(R)$ in the domain and restricted to $Cen(MC_n(R))$ in the codomain. So first we should see that $i_n(Cen(R)) \subseteq Cen(MC_n(R))$ to verify that $f$ is well defined, then let $\hat{i}^n_k(A) \in MC_n(R)$ with $A \in M_{n^k}(R)$ and $a \in Cen(R)$. Therefore $i_n(a) \hat{i}^n_k(A) = \hat{i}^n_k(i_n(a)A)$ and $(i_n(a)A)_{ij} = a A_{ij} = A_{ij} a = (A i_n(a))_{ij}$ for $1 \leq i, j \leq n$ which means that $i_n(a) \hat{i}^n_k(A) = \hat{i}^n_k(A) i_n(a)$, so $i_n(a) \in Cen(MC_n(R))$. Also we note as the function $i_n$ is injective the function $f$ is injective, so consider $\hat{i}^n_k(A) \in Cen(MC_n(R))$ with $A \in M_{n^k}(R)$ then $A \in Cen(M_{n^k}(R))$ and $\alpha[n]^k_1(a) = A$ for some $a \in Cen(R)$ then $f(a) = i_n(a) = \hat{i}^n_k(A)$ as desired.

\begin{proposition}
Let $n$ be a positive integer. Then $U(MC_n(R)) = \bigcup_{k=0}^\infty \hat{i}^n_k(U(M_{n^k}(R)))$.
\end{proposition}

\textbf{Proof.} As the multiplication behaves in a local way, an element has inverse if and only if it has inverse in its original matrix ring.

\begin{proposition}
Let $n$ be a positive integer. Then $|MC_n(R)| = \max\{|R|, \aleph_0\}$.
\end{proposition}

\textbf{Proof.} The cardinality of $MC_n(R)$ is the supremum of the cardinalities of $M_{n^k}(R)$ which is $\max\{|R|, \aleph_0\}$.

\section{Functors}

Let $n$ be a positive integer and let $R$ be a ring we may consider the following family of functors $\delta : R\text{-Mod} \to M_n(R)\text{-Mod}$ given by $\delta(M) = M^n$ for $M$ a $R$-module and $M_n(R)$ action given by $(A \phi)(k) = \sum_{i=1}^n A_{ik} \phi(k)$ for $A \in M_n(R)$ and $\phi \in M^n$, and for any $f : M \to N$ we define $\delta(f) = f^n : M^n \to N^n$. Now we define $\delta_0 = \delta$ and $\delta_{k+1} = \delta \circ \delta_k$, so we obtain a direct system of functors and we put $\Delta = \varinjlim \delta_k : R\text{-Mod} \to MC_n(R)\text{-Mod}$. Note that $\delta_k$ is the usual functor that give us the Morita equivalence between $R\text{-Mod}$ and $M_{n^k}(R)\text{-Mod}$ so the functor $\Delta$ is in a way an attempt to get and equivalence between $R\text{-Mod}$ and $MC_n(R)\text{-Mod}$, almost never $\Delta$ is an equivalence but there is information that we may obtain from $\Delta$.

\begin{proposition}
Let $n$ be a positive integer. Then $\Delta : R\text{-Mod} \to MC_n(R)\text{-Mod}$ is a functor.
\end{proposition}

For $M$ a left $R$-module an element of $\Delta(M)$ is given by $i^k_M(m)$ where $k$ is a natural and $m \in M^{n^k}$.

\begin{proposition}
Let $n$ be a positive integer, let $R$ be a ring and let $M$ be a left $R$-module. Then $i^k_M(M^{n^k})$ is a $MC_n(R)$-generator of $\Delta(M)$.
\end{proposition}

\begin{proposition}
Let $n$ be a positive integer, let $R$ be a ring and let $S$ be a left simple $R$-module. Then $\Delta(S)$ is a left simple $MC_n(R)$-module.
\end{proposition}

\textbf{Proof.} Let $i^k_S(m)$ be a no zero element in $\Delta(S)$ with $k$ a natural and $m \in S^{n^k}$ then $i^k_S(m)$ $MC_n(R)$-generates $i^k_M(S^{n^k})$, since it $M_{n^k}$-generates it. Finally $i^k_M(S^{n^k})$ is contained in the $MC_n(R)$-submodule $MC_n(R)$-generated by $i^k_S(m)$ the proposition follows.

\begin{proposition}
Let $n$ be a positive integer, let $R$ be a ring and let $M$ be a left faithful $R$-module. Then $\Delta(M)$ is a left faithful $MC_n(R)$-module.
\end{proposition}

\textbf{Proof.} Let $i^n_k(A)$ be an element in $\text{Ann}_{MC_n(R)}(\Delta(M))$ with $k$ a natural and $A \in M_{n^k}(R)$, so $A \in \text{Ann}_{M_{n^k}(R)}(i^k_M(M^{n^k}))$ as $M^{n^k}$ is $M_{n^k}(R)$ faithful then $A = 0$ and $\text{Ann}_{MC_n(R)}(\Delta(M)) = 0$.

\begin{corollary}
Let $n$ be a positive integer and let $R$ be a ring. If $R$ is a left primitive ring then $MC_n(R)$ is a left primitive ring.
\end{corollary}

\begin{proposition}
Let $n$ be a positive integer and let $R$ be a ring. Then $\Delta$ commutes with coproducts.
\end{proposition}

\begin{proposition}
Let $n$ be a positive integer and let $R$ be a ring. Then $\Delta(R)$ is a free left $R$-module.
\end{proposition}

\begin{corollary}
Let $n$ be a positive integer and let $R$ be a ring. If $P$ is projective left $R$-module, then $\Delta(P)$ is a projective left $MC_n(R)$-module.
\end{corollary}

There is another functor from $R$-Mod to $MC_n(R)$-Mod, first consider the ring morphism $i_R : R \to MC_n(R)$ the canonical inclusion, so any left $MC_n(R)$-module has structure of left $R$-module, in particular and by symmetry, $MC_n(R)$ has structure of $R$-$R$-bimodule, and so we may define the new functor as $\Phi(M) = MC_n(R) \otimes_R M$ for any left $R$-module $M$ and $\Phi(f) = 1_{MC_n(R)} \otimes f$ for any left $R$-morphism $f$. There is naive way to compare the two functors.

\begin{proposition}
Let $n$ be a positive integer and let $R$ be a ring. Then $MC_n(R)$ is a free left $R$-module.
\end{proposition}

\begin{corollary}
Let $n$ be a positive integer, let $R$ be a ring and let $M$ be a left $R$-module. If $\text{pd}_R(M) = k$ then $\text{pd}_{MC_n(R)}(\Phi(M)) \leq k$.
\end{corollary}

\begin{proposition}
Let $n$ be a positive integer and let $R$ be a ring. There is a natural transformation $\Pi : \Phi \to \Delta$.
\end{proposition}

\textbf{Proof.} Let $k$ be a natural and let $M$ be a left $R$-module. Then we define a $M_{n^k}(R)$-morphism $\Pi^k_M : M_{n^k}(R) \otimes_R M \to \delta_k$ as $\Pi^k_M(A \otimes m) = A \delta_k(m)$ for any $A \in M_{n^k}(R)$ and any $m \in M$. In fact $\Pi^k$ is a natural transformation from $R$-Mod to $M_{n^k}(R)$-Mod, so we may take the direct limit of the direct system $\{\Pi^k\}_{k \in \mathbb{N}}$ and that is the natural transformation $\Pi$ we are looking for.

\section{Ideals of the Matricial Closure}

The direct system created for a ring $R$ serves as direct system for any left (right, bilateral) ideal $I$ of $R$, at least in the category of abelian groups, so we have that $MC_n(I)$ is a subgroup of $MC_n(R)$. Moreover we obtain:

\begin{proposition}
Let $n$ be a positive integer, let $R$ be a ring and let $I$ be a left (right, bilateral) ideal of $R$. Then $MC_n(I)$ is a left (right, bilateral) ideal of $MC_n(R)$.
\end{proposition}

\textbf{Proof.} First we name $w_k : M_{n^k}(I) \to MC_n(I)$ the canonical inclusions of the matricial closure of a left (right, bilateral) ideal with $k$ a natural and $\beta^m_k : M_{n^k}(I) \to M_{n^m}(I)$ the connection morphisms with $k \leq m$ naturals. Let $w_k(X) \in MC_n(I)$ and $i^n_m(A)$ with $X \in M_{n^k}(I)$ and $A \in M_{n^k}(R)$, we begin with the case when $k \leq m$, in this case we put $i^n_m(A) w_k(X) := w_m(A \beta^m_k(X))$ and when $k > m$ we put $i^n_m(A) w_k(X) := w_k(\alpha^k_m(A) X)$, with this left action defined $MC_n(I)$ becomes a left ideal of $MC_n(R)$. In the same manner for right and bilateral ideal. $\blacksquare$

The next proposition tell us that the matricial closure behaves like the matrix ring respect to the bilateral ideals.

\begin{proposition}
Let $n$ be a positive integer and let $R$ be a ring. Then there is a lattice isomorphism between the bilateral ideals of $R$ and the bilateral ideals of $MC_n(R)$.
\end{proposition}

\textbf{Proof.} We take an ideal $I$ of $MC_n(R)$ and define the next subset of $R$, $J = \{A_{11} \in R \mid i^n_k(A) \in I, A \in M_{n^k}(R)\}$, it is easy to see that $J$ is an ideal of $R$. Next we call $\{e^k_{ij}\}_{1 \leq i, j \leq k}$ the canonical basis of $M_k(R)$, and for $A \in M_k(R)$ we remark that $(A e^k_{rs})_{ij} = A_{ir} \delta_{sj}$ and $(e^k_{rs} A)_{ij} = \delta_{ir} A_{sj}$ with $1 \leq i, j, r, s \leq k$ and $\delta$ stands for the Kronecker's delta. We wish to prove that $MC_n(J) = I$. Let $w_k(X) \in MC_n(J)$ with $X \in M_{n^k}(J)$ and we write
\[
X = \sum_{i=1}^{n^k} \sum_{j=1}^{n^k} X_{ij} e^{n^k}_{ij}
\]
For all $1 \leq i, j \leq n^k$ there exists a matrix $Y^{ij} \in M_{n^{k_{ij}}}(R)$ with $k_{ij}$ a natural such that $w_{k_{ij}}(Y^{ij}) \in I$ and $Y^{ij}_{11} = X_{ij}$, cause there is a finite number of $k_{ij}$ with a lose of generality we may think via the connection morphism that $k = k_{ij}$, making $k$ bigger if we need it. Now as a simple fact
\[
e^{n^k}_{i1} Y^{ij} e^{n^k}_{1j} = X_{ij} e^{n^k}_{ij}
\]
for all $1 \leq i, j \leq n^k$. Applying $w_k$ and summing, we get
\[
\sum_{i=1}^{n^k} \sum_{j=1}^{n^k} w_k(e^{n^k}_{i1}) w_k(Y^{ij}) w_k(e^{n^k}_{1j}) = w_k\left(\sum_{i=1}^{n^k} \sum_{j=1}^{n^k} X_{ij} e^{n^k}_{ij}\right) = w_k(X)
\]
As $I$ is an ideal of $MC_n(R)$, then $w_k(X) \in I$. Now let $w_k(X) \in I$ with $X \in M_{n^k}(R)$ and take
\[
X = \sum_{i=1}^{n^k} \sum_{j=1}^{n^k} X_{ij} e^{n^k}_{ij}
\]
and if we observe that $X_{ij} = (e^{n^k}_{1i} X e^{n^k}_{j1})_{11}$ so $X_{ij} \in J$ for any $1 \leq i, j \leq n^k$, and so we have that the assignation $MC_n$ from the ideals of $R$ to the ideals of $MC_n(R)$ is onto, is easy to observe that is also into and monotone, so it is a lattice isomorphism. $\blacksquare$

\begin{corollary}
For any positive $n$ integer and ring $R$, if $R$ is a simple ring then $MC_n(R)$ is a simple ring.
\end{corollary}

\begin{proposition}
Let $n$ be a positive integer, let $R$ be a ring and let $\{I_k\}_{k=0}^\infty$ be a family such that $I_k$ is a left (right, bilateral) ideal of $M_{n^k}(R)$ with $w_k(I_k) \subseteq w_{k+1}(I_{k+1})$. Then $I = \bigcup_{k=0}^\infty w_k(I_k)$ is a left (right, bilateral) ideal of $MC_n(R)$. Moreover, if $I_k$ is a maximal left (right, bilateral) ideal of $M_{n^k}(R)$ then $I$ is maximal left (right, bilateral) ideal of $MC_n(R)$.
\end{proposition}

\textbf{Proof.} The proof follows from the fact that $i^n_k(R)$ is a left (right, bilateral) $MC_n(R)$-generator set of $MC_n(R)$ for all $k$ natural. $\blacksquare$

Now let $K$ be a field and $k, m$ be naturals with $0 \leq m < 2^k$, and we define $I^k_m = \{A \in M_{2^k}(K) \mid A_{im} = 0 \text{ for all } i = 1, \dots, 2^k\}$, note that $I^k_m$ is a maximal left ideal of $M_{2^k}(K)$. We proceed by consider the binary expansion of $m = (m_0, \dots, m_{k-1})$ with $m_i \in \{0, 1\}$ for all $i = 0, \dots, k-1$. Observe that for $m, m'$ naturals such that $0 \leq m < 2^k$ and $0 \leq m' < 2^{k+1}$ with binary expansion $m = (m_0, \dots, m_{k-1})$ and $m' = (m'_0, \dots, m'_k)$, $\beta^{k+1}_k(I^k_m) \subset I^{k+1}_{m'}$ if and only if $m'_i = m_i$ for all $i = 0, \dots, k-1$, which means that if we take a function $\omega : \mathbb{N} \to \{0, 1\}$ this defines a family of left maximal ideals with the condition above so $I_\omega = \bigcup_{k=0}^\infty w_k(I^k_{(\omega(0), \dots, \omega(k-1))})$ is a maximal left ideal of $MC_2(K)$. It is standard knowledge that the quantity of this functions is $2^{\aleph_0}$, then we obtained $2^{\aleph_0}$ maximal left ideals. Consider the next proposition.

\begin{proposition}
Let $R$ be a ring and let $S, S'$ be simple left $R$-modules such that $S \cong R/I$ and $S' \cong R/I'$. Then $S \cong S'$ if and only if there is an $a \in R$ such that $I' = (I : a)$.
\end{proposition}

If we consider an enumerable field such as $\mathbb{Z}_2$, $MC_2(K)$ has $\aleph_0$ elements, so for any maximal left ideal $I$ of $MC_2(K)$ there is $\aleph_0$ left maximal ideals of $MC_2(K)$ that induce the same left simple $MC_n(K)$-module, as there is at least $2^{\aleph_0}$ maximal left ideals of $MC_2(K)$ so there is at least $2^{\aleph_0}$ simple left $MC_2(K)$-modules non-isomorphic.
\section{Chain Conditions}

\subsection{Descending Chain Condition}

All left artinian rings satisfy the descending chain condition on left direct summands. Let $n$ be a positive integer and let $R$ be a ring we denote $\{e^{n}_{ij}\}_{1 \leq i, j \leq n}$ the canonical basis of $M_n(R)$, we may build the next strictly descending chain of left direct summands $I_k = MC_n(R) i^n_k(e^{n^k}_{11})$ for any natural $k$. First note as $i^n_k(e^{n^k}_{11})$ is an idempotent, then $I_k$ is a left direct summand of $MC_n(R)$, also as $i^n_k(e^{n^k}_{11})$ divides by the left to $i^n_k(e^{n^{k+1}}_{11})$ the chain is descending, at last $\bigcap_{m \in \mathbb{N}} I_m = 0$ and this with the fact that $i^n_k(e^{n^k}_{11}) \in I_k$ implying $I_k \neq 0$ means that $MC_n(R)$ does not satisfy the descending condition on left direct summands, which means $MC_n(R)$ is not artinian, neither. So we have.

\begin{proposition}
Let $n$ a positive integer and let $R$ be a no zero ring. Then $MC_n(R)$ does not satisfy descending chain condition on left direct summands.
\end{proposition}

\begin{corollary}
Let $n$ a positive integer and let $R$ be a no zero ring. Then $MC_n(R)$ is not artinian.
\end{corollary}

\begin{corollary}
Let $n$ a positive integer and let $R$ be a no zero ring. Then $MC_n(R)$ is not semisimple artinian.
\end{corollary}

\begin{corollary}
Let $n$ a positive integer and let $R$ be a no zero ring. Then $MC_n(R)$ does not have left finite uniform dimension.
\end{corollary}

As if a ring has finite left uniform dimension then it has left Krull dimension, the next corollary follows:

\begin{corollary}
Let $n$ a positive integer and let $R$ be a no zero ring. Then $MC_n(R)$ does not have left Krull dimension.
\end{corollary}

\begin{corollary}
Let $n$ a positive integer and let $R$ be a no zero simple ring. Then $MC_n(R)$ is not a left full linear ring.
\end{corollary}

\textbf{Proof.} As $\text{Soc}(MC_n(R)_{MC_n(R)})$ is a bilateral ideal and $MC_n(R)$ is simple then $\text{Soc}(MC_n(R)_{MC_n(R)})$ is $MC_n(R)$ or $0$, as the first option will imply that $MC_n(R)$ is semisimple artinian which never happens then it should be the second one.

\subsection{Ascending Chain Condition}

First as we observe $MC_n(R)$ never is semisimple artinian, unless $R = 0$, as $MC_n(R)$ preserves the property of being Von Neumann regular, $MC_n$ does not preserve the property of being left noetherian, as any left noetherian and von Neumann regular is semisimple artinian. As example for any field $K$, $MC_n(K)$ is not left noetherian but $K$ is. Moreover the descending chain created below give us a strictly ascending chain of direct left summands. So we get:

\begin{corollary}
Let $n$ a positive integer and let $R$ be a no zero ring. Then $MC_n(R)$ does not satisfy ascending chain condition on left direct summands.
\end{corollary}

\begin{corollary}
Let $n$ a positive integer and let $R$ be a no zero ring. Then $MC_n(R)$ is not noetherian.
\end{corollary}

\begin{corollary}
Let $n$ a positive integer and let $R$ be a no zero ring. Then $MC_n(R)$ is not quasi Frobenius.
\end{corollary}

\section{Jacobson Radical and Matrix Closure}

We are going prove some results in the spirit of radical theory of rings to demonstrate in a easy way that certain properties are preserved by $n$-matricial closure.

\begin{proposition}
For any positive integer $n$ and any ring $R$, $MC_n(J(R)) = J(MC_n(R))$.
\end{proposition}

\textbf{Proof.} For the characterization mentioned in the preliminaries, $J(MC_n(R))$ is the unique left ideal with all its elements left quasiregular and maximal with this property. So note that $MC_n(J(R))$ is an ideal of $MC_n(R)$ in particular a left ideal, also any element is of the form $i^n_k(A)$ with $A \in M_{n^k}(J(R))$ and $k$ a natural, as $M_{n^k}(J(R)) = J(M_{n^k}(R))$, there $1 - i^n_k(A)$ is left invertible then $i^n_k(A)$ is left quasiregular, which means that all elements of $MC_n(J(R))$ are left quasiregular, therefore $MC_n(J(R)) \subseteq J(MC_n(R))$. At last notice that a left quasiregular element in $MC_n(R)$ which is in $i^n_k(M_{n^k}(R))$ is left quasiregular in $i^n_k(M_{n^k}(R))$, so $i^n_k(M_{n^k}(R)) \cap J(MC_n(R)) = i^n_k(M_{n^k}(J(R)))$ and we got the proposition.

\begin{corollary}
Let $n$ a positive integer and let $R$ be a no zero ring. If $R$ is semisimple then $MC_n(R)$ is semisimple.
\end{corollary}

\begin{corollary}
Let $n$ a positive integer and let $R$ be a no zero ring. Then $MC_n(R)$ is not semiperfect.
\end{corollary}

\begin{corollary}
Let $n$ a positive integer and let $R$ be a no zero ring. Then $MC_n(R)$ is not perfect.
\end{corollary}

\section{Certain Classes of Rings and Matrix Closure}

\subsection{Von Neumann Regular Rings}

\begin{proposition}
Let $n$ a positive integer and let $R$ be a no zero ring. If $R$ is von Neumann regular ring then $MC_n(R)$ is von Neumann regular ring.
\end{proposition}

\textbf{Proof.} Let $i^n_k(A)$ be an element in $MC_n(R)$ with $A \in M_{n^k}(R)$, as $M_{n^k}(R)$ is a von Neumann regular ring then there is $X \in M_{n^k}(R)$ with $A = A X A$, so $i^n_k(A) = i^n_k(A) i^n_k(X) i^n_k(A)$ as desired.

\subsection{V-rings}

As Herbera proves in her paper [3], for any von Neumann regular ring $R$, if it has dimension less than $2^{\aleph_0}$ over its center, then there is no left injective simple modules over $R$. We know that for $K$ a field and $n$ a positive integer, if we put $R = MC_n(K)$ its center is isomorphic to $K$ and also that the dimension of $R$ is exactly $\aleph_0$. So the next result follows:

\begin{proposition}
Let $n$ a positive integer and let $K$ be a field. Then $MC_n(K)$ is not a V-ring.
\end{proposition}

\subsection{Semiprime Rings}

\begin{proposition}
Let $n$ be a positive integer. If $R$ is a semiprime ring then $MC_n(R)$ is a semiprime ring.
\end{proposition}

\textbf{Proof.} First we recall that if $R$ is semiprime then $M_m(R)$ is semiprime for any positive integer $m$. Now let $i^n_k(A), i^n_k(X) \in MC_n(R)$ with $i^n_k(A) i^n_k(X) i^n_k(A) = 0$ and without lose of generality with $A, X \in M_{n^k}(R)$, so if $i^n_k(A X A) = i^n_k(A) i^n_k(X) i^n_k(A) = 0$ then $A X A = 0$, since $M_{n^k}(R)$ is semiprime then $A = 0$, meaning that $i^n_k(A) = 0$.

\subsection{Prime Rings}

\begin{proposition}
Let $n$ be a positive integer. If $R$ is a prime ring then $MC_n(R)$ is a prime ring.
\end{proposition}

\textbf{Proof.} As above we recall that if $R$ is prime then $M_m(R)$ is prime for any positive integer $m$. Now let $\hat{i}^n_k(A), \hat{i}^n_k(X), \hat{i}^n_k(B) \in MC_n(R)$ with $\hat{i}^n_k(A) \hat{i}^n_k(X) \hat{i}^n_k(B) = 0$ and without lose of generality with $A, X, B \in M_{n^k}(R)$, so if $\hat{i}^n_k(A X B) = \hat{i}^n_k(A) \hat{i}^n_k(X) \hat{i}^n_k(B) = 0$ then $A X B = 0$, since $M_{n^k}(R)$ is prime then $A = 0$ or $B = 0$, it follows that $\hat{i}^n_k(A) = 0$ or $\hat{i}^n_k(B) = 0$.
\section{K-Theory}

\subsection{Invariant Basis Number}

Let $n$ and $k$ be natural numbers and let $K$ be a field, we note that the group morphism $K_0(\alpha[n]^{k+1}_k) : K_0(M_{n^k}(K)) \to K_0(M_{n^{k+1}}(K))$ is the multiplication by $n$, then from the fact that $K_0$ commute with direct limits and $K_0(M_{n^m}(K)) = \mathbb{Z}$ for any $m$ natural, we obtain that $K_0(MC_n(K))$ is the direct limit of the of the endomorphism multiplication by $n$ of $\mathbb{Z}$, and the limit is $\mathbb{Z}[\frac{1}{n}]$. So from this we get that $MC_n(K)$ have the invariant basis number. In general is unknown if the functor preserves the invariant basis number.

\section{Bibliography}

\begin{enumerate}
\item P. M. Cohn, Free rings and their relations. Second edition. London Mathematical Society Monographs, 19. Academic Press, Inc. [Harcourt Brace Jovanovich, Publishers], London, 1985.
\item K. R. Goodearl, Von Neumann regular rings. Second edition. Robert E. Krieger Publishing Co., Inc., Malabar, FL, 1991.
\item D. Herbera, Simple modules over small rings. Rings, modules and representations, 189205, Contemp. Math., 480, Amer. Math. Soc., Providence, RI, 2009.
\item T. Y. Lam, Lectures on modules and rings. Graduate Texts in Mathematics, 189. Springer-Verlag, New York, 1999.
\item J. C. McConell and J. C. Robson, Noncommutative noetherian rings. With the cooperation of L. W. Small. Revised edition. Graduate Studies in Mathematics, 30. American Mathematical Society, Providence, RI, 2001.
\item C. Nastessecu, F. Van Oystaeyen, Dimensions of ring theory. Mathematics and its Applications, 36. D. Reidel Publishing Co., Dordrecht, 1987.
\item J. Rosenberg, Algebraic K-theory and its applications. Graduate Texts in Mathematics, 147. Springer-Verlag, New York, 1994.
\item R. Wisbauer, Foundations of Module and Ring Theory. Gordon and Breach, Reading 1991.
\end{enumerate}

\end{document}